\documentclass[12pt]{article}
\usepackage{amsmath,amsthm}
\usepackage{amssymb,latexsym}
\usepackage{amsfonts}

\newtheorem{thm}{Theorem}[section]
\newtheorem{cor}{Corollary}[section]

\theoremstyle{definition}

\begin{document}

\title{Homology of contact CR-submanifolds}
\author{Fulya \c{S}ahin and Bayram \c{S}ahin \\
\vspace{6pt}{\small{\it Ege University, Department of Mathematics, 35100, Izmir, Turkey}}\\
\vspace{6pt}{\small{\it fulya.sahin@ege.edu.tr}, {\it bayram.sahin@ege.edu.tr}}}
\date{}
\maketitle

\begin{abstract}
\baselineskip=16pt

In this paper, homology of a contact CR-submanifold of a real hypersurface, which has naturally almost contact metric structure induced from  the complex Euclidean space $\mathbb{C}^{m}$, is examined. More precisely, nonexistence of stable integral currents on a compact contact CR-submanifold of  real hypersurface of canonical complex space form $\mathbb{C}^m$ is investigated and  vanishing
theorems concerning the homology groups are obtained.

 \end{abstract}
\noindent{\small {\bf Mathematics Subject Classifications (2010).} 55N25, 57R60, 53C65, 53B25.}\\
\noindent{\small {\bf Key words.} Almost contact manifold, real hypersurface, sphere theorem, CR-warped submanifold, stable current, homology group.} \\

\section{Introduction}

Recently, geometric and topological methods have important roles in learning neural networks, pattern recognition and signal processing \cite{BBCSV}, \cite{UY}. The main idea in these areas is to consider a flat or curved network parameter space endowed with  a suitable geometric structure in the network learning algorithm. The non-Euclidean nature of  data implies that manifold theory may be considered as new tool to use in these new research areas. One of the new methods in machine learning is depicted  data as a submanifold (with curvature) of Euclidean space. Thus submanifold theory finds  a place in machine learning and data mining, which are prominent research areas.
Let $\tilde{M}$ be an $m$-dimensional compact Riemannian manifold with Riemannian metric
$\tilde{g}$ and the Levi-Civita connection $\tilde{\nabla}$. Denote by $(S, \varrho)$ the oriented, $p$-rectifiable set in
$\tilde{M}$. The set of rectifiable $p$-currents is
$$\mathcal{R}_p(\tilde{M})=\{\mathcal{S}:\sum^{\infty}_{n=1}n\mathcal{S}_n; \mathcal{S}_n=( S_n, \varrho_n), \tilde{M}(\mathcal{S})=\sum^{\infty}_{n=1}n\mathcal{H}^p(S_n)<\infty\}.$$
$\mathcal{S}\in  \mathcal{R}_p(\tilde{M})$ is called an integral $p$-current if $\mathcal{S}$ and $\partial \mathcal{S}$ are both rectifiable currents, for details see:\cite{GM}. After Federer
and Fleming \cite{FF} showed that any non-trivial integral homology
class in $H_p(M,\mathbb{Z})$ corresponds to  a stable current, Lawson and Simons \cite{LS} obtained that there are no stable
integral currents in the sphere $S^n$, and there is no integral
current in a submanifold $M^m$ of $S^n$ when the second fundamental
form of $M^m$ satisfies a pinching condition. The result on
submanifolds of $S^n$ has been extended to Euclidean space \cite{Leung}, \cite{Xin}, to $S^{n_1}\times
S^{n_2}$  \cite{Zhang} and to higher codimensional submanifolds
in the Euclidean spaces \cite{LZ}. In \cite{CS}, the authors showed that when the second fundamental form of a compact  submanifold $N^n$ of a hypersurface $M^{m}$ in $\mathbb{R}^{m+1}$ or  a submanifold $M^m$ in $\mathbb{R}^{m+2}$ satisfies some
conditions, there is no stable integral $p-$ current in $N^n$ and $H_p(N, \mathbb{Z})= H_{n-p}(N, \mathbb{Z}) =0$,  where
$H_i(N,\mathbb{Z})$ is the $i-$ th homology group of $N$ with
integer coefficients.

CR-submanifolds of a K\"{a}hler manifold were introduced by Bejancu \cite{Bejancu}. The study of such submanifolds is still an active research area, see:\cite{Vilcu2020}, \cite{Vilcu} for instances.  On the other hand, Chen defined CR-warped product submanifolds of
almost Hermitian manifolds in his papers \cite{CCR1} and \cite{CCR2}. Later such
submanifolds have been studied by many authors, see  book
\cite{ChenB} and references therein. In almost contact manifolds, Hasegawa and Mihai \cite{HM1}  introduced and studied contact CR-warped submanifolds and  obtained a sharp
inequality for the squared norm of the second fundamental form (an extrinsic invariant) in
terms of the warping function for contact CR-warped products isometrically immersed in
Sasakian manifolds. Moreover Mihai  \cite{M} obtained a classification of contact
CR-warped product submanifolds in spheres, which satisfy the equality case, see also  \cite{M2, MM, Munteanu}.

In \cite{SahinF}, the first author considered CR-warped product submanifolds
in even-dimensional Euclidean spaces and obtained a sufficient condition for certain homology groups to be zero. In this article, a more general case is taken into account in calculating homology groups of CR-submanifolds immersed in an arbitrarily curved space. In this direction,  first the non-existence of integral currents on a CR-warped submanifold of the real hypersurface of complex Euclidean space was found in terms of the warping function. Second, the non-existence of integral currents on a  general contact CR-submanifold of the real hypersurface of complex Euclidean space was obtained by using  mixed scalar curvature introduced in \cite{R1}, \cite{RZ} and related notions.

\section{Preliminaries}
In this section, we will review basic materials from \cite{Yano-Kon} for later sections.
An almost contact manifold $M$ is a manifold endowed with an $(1,1)$ tensor field $\phi$, a 1- form $\eta$ and a vector field $\xi$ such that
\begin{equation}
\phi^{2}X=-X+\eta(X)\xi,  \eta(\xi)=1, \phi\xi=0 \label{eq:2.1}
\end{equation}
for $X\in \Gamma(TM)$, where $\Gamma(TM)$ denotes the module of vector fields on $M$. An almost contact metric structure $(\phi,\xi,\eta)$ is called Sasakian structure if
\begin{equation}
(\nabla_{X}\phi)Y=g(X,Y)\xi-\eta(Y)X \nonumber
\end{equation}
for $Y\in \Gamma(TM)$, where $\nabla$ is the Levi-Civita connection on $M$. A K\"{a}hler manifold $\bar{M}$ is a differentiable Riemannian manifold endowed with an almost complex structure $J$, and a Riemannian metric $g$ such that
\begin{eqnarray}
J^{2}X=-X,  \nabla J=0,  g(JX,JY)=g(X,Y) \label{eq:2.2}
\end{eqnarray}
for $X, Y\in \Gamma(T\bar{M})$. A complex space form is a simply connected complete K\"{a}hler manifold of constant holomorphic sectional curvature.

Let $M$ be an orientable real hypersurface of the complex space form $\mathbb{C}^{m}$, which is a complex space form with constant holomorphic sectional curvature $0$, with the unit normal vector field ${\zeta}$. We denote the induced metric on $M$ by $g$ again. Then there is an almost contact metric structure $(\phi,\xi,g)$ defined on $M$, whose $\xi$ is the characteristic unit vector field defined by $\xi=-J{\zeta}$, $\eta$ is a smooth form which is dual to $\xi$ that is $\eta(X)=g(X,\xi)$, $X\in \Gamma(TM)$ and $\phi$ is a tensor field of type $(1,1)$ satisfying (\ref{eq:2.1}) and
$$
g(\phi X,\phi Y)=g(X,Y)-\eta(X)\eta(Y).
$$
Let $A_{\zeta} $ be the shape operator of $M$ with respect to $\zeta$. Then we also have the following identities
\begin{eqnarray}
R(X,Y)Z=g(AY,Z)AX-g(AX,Z)AY\label{eq:2.6},\\
JX=\phi X+\eta(X)\xi\label{eq:2.8}
\end{eqnarray}
where we have put $A=A_{\zeta}$.
Let $N$ be an $n+1$- dimensional submanifold, tangent to the structure vector field $\xi$ of a real hypersurface $M^{2m-1}$ of $\mathbb{C}^{m}$. The Gauss equation for a submanifold $N$ is given by
\begin{eqnarray}
  g(R(X,Y)Z,W)&=&g(R'(X,Y)Z,W)-g(h'(X,W),h'(Y,Z))\nonumber\\
  &&+g(h'(Y,W),h'(X,Z))\label{eq:2.15}
\end{eqnarray}
 for $X,Y,Z,W \in \Gamma(TN)$, here $h'$ is
 the second fundamental form of $N$, $R'$ is the curvarture tensor field of $N$ and $R$ is the curvature tensor field of $M$.
From (\ref{eq:2.1}), (\ref{eq:2.2}) and (\ref{eq:2.8}), it follows that a real hypersurface of a K\"{a}hler manifold is Sasakian manifold if and only if $A$ satisfies
\begin{equation}
A=-I+\beta \eta\otimes \xi \nonumber
\end{equation}
where $\beta=\eta(A\xi)+1.$ For research papers on the real hypersurfaces of a Kaehler manifold being a Sasakian manifold, see: \cite{Okumura} and \cite{Tashiro}, see also  \cite{Carriazo} for generalized Sasakian space form.

We now suppose $N$ is a contact CR-submanifold of $M^{2m-1}$, then there are two distrributions $D$ and $D^{\perp}$ on $N$ such that
\begin{eqnarray}
TM=D\oplus D^{\perp}\oplus {\xi}, \phi(D)=D, \phi(D^{\perp})\subset T^{\perp}M.\nonumber
\end{eqnarray}
We note that totally real submanifolds and invariant submanifolds are contact CR-submanifolds with $D=\{0\}$ and $D^{\perp}=\{0\}$, respectively. For examples of such submanifolds, see: \cite{Yano-Kon2}. Unless otherwise stated, we will consider proper contact CR-submanifolds throughout this paper. This means that $dim(D)=p\neq 0$ and $dim(D^{\perp})=q\neq 0$.

 \section{Homology of Contact CR-Warped Submanifolds of a Sasakian Hypersurface}
    CR-warped product submanifolds of a K\"{a}hler manifold were first introduced by Chen \cite{ChenB, CCR1, CCR2} and then such submanifolds were studied by Hasegawa and Mihai\cite{HM1}, Mihai \cite{M} and Munteanu \cite{Munteanu} in the Sasakian setting. We first recall basic notions for warped product manifolds from \cite{BO} and \cite{ChenB}. Let $(B, g_1)$ and $(F, g_2) $ be two Riemannian manifolds, $f:B
\rightarrow (0,\infty) $ and $\pi: B \times F \rightarrow B$,
$\eta:B \times F \rightarrow F$ the projection maps given by $\pi
(p,q)=p$ and $\eta (p,q)=q$ for every $(p,q) \in B \times F$.  The
warped product \cite{BO} $M=B\times_{f} F$ is the manifold $B
\times F$ equipped with the Riemannian structure such that
$$ g(X,Y)=g_1(\pi_* X,\pi_* Y)+(f o \pi)^2
g_2(\eta_* X,\eta_* Y)$$ for every $X$ and $Y$ of $M$, where $*$
denotes the tangent map. The function $f$ is called the warping
function of the warped product manifold. In particular, if the
warping function is constant, then the warped product manifold $M$
is said to be trivial. Let $ X, Y$ be vector fields on $B$ and $V,W$
vector fields on $F$, then from Lemma~7.3 of \cite{BO}, we have
\begin{equation}
\nabla_X V=\nabla_V X=(\frac{Xf}{f})V \nonumber
\end{equation}
where $\nabla$ is the Levi-Civita connection on $M$.

    In \cite{HM1}, Hasegawa and Mihai considered warped products $M=M_{1}\times_fM_2$ which are contact CR-submanifolds of a Sasakian manifold $\bar{M}$. They showed  that there do not exist warped product submanifolds $M=M_{1}\times_fM_2$ of a Sasakian manifold $\bar{M}$ such that $M_1$ is an anti-invariant submanifold tangent to $\xi$ and $M_2$ an invariant submanifold of $\bar{M}$. Therefore they considered contact  CR-warped product  submanifold of the  form $M_{1}\times_fM_2$ of $\bar{M}$ such that $M_1$ is an invariant submanifold to $\bar{M}$ and $M_2$ an anti-invariant submanifold of $\bar{M}$.

One of the generalized version of the Lawson-Simons's result is obtained as follows.
\begin{thm}\label{Cheng-Shiohama}(Cheng-Shiohama)\cite{CS} Let $N$ be an $n$-dimensional compact submanifold of $M^m$ which is an $m$-dimensional hypersurface of an $(m+1)$-dimensional Euclidean space $\mathbb{R}^{m+1}$. If
\begin{eqnarray*}
\sum_{i,\alpha}[2\|h'(e_i,e_\alpha)\|^2-\langle h'(e_i,e_i),h(e_\alpha,e_\alpha)\rangle]<2p(p-n)\lambda^2
\end{eqnarray*}
is satisfied for every $x\in N$ and any orthonormal basis $\{e_i,e_\alpha\}$ of $T_xN$, $i=1,...,p$, $\alpha=p+1,...,n.$ Here $h'$ is the second fundamental form of $N$ in $M^m$, and $\lambda^2$ is the maximum at $x$ of the square of principal curvatures of $M^m$. Then there is no integral $p$-current in $N$ and $H_p(N,\mathbb{Z})=H_{n-p}(N,\mathbb{Z})=0$, where $H_p(N,\mathbb{Z})$ is the $p$-th singular homology group with integer coefficients.
\end{thm}

In this section, we are going to find suitable conditions for a contact CR-warped product such that $(2p+1)th$ homology group of $N$ is zero and as a result of this, we obtain a sphere theorem for such submanifolds. We note that the Laplacian of a function $f$ is denoted by $ \bigtriangleup f$.
     \begin{thm}
    Let $N^{n}$ be a compact contact CR-warped product submanifold $M=M_1\times _fM_2$ of a Sasakian hypersurface of $\mathbb{C}^{m}$ such that $M_1$ is a $(2p+1)$-dimensional invariant submanifold tangent to $\xi$ and $M_2$ a $q$-dimensional $\mathcal{C}$- totally real submanifold of $M$. If
    \begin{eqnarray*}
    \bigtriangleup f&<&\frac{f}{q}\{2(2p+1)(2p+1-n)\lambda (x)^2-q(2p-\eta(A\xi))+\|\nabla lnf\|^2+1)\}
    \end{eqnarray*}
    is satisfied for every $x\in N$, where $\lambda (x)^2$ is the maximum at $x$ of the square of principal curvature of $M^{2m-1}$, then  there is no stable integral $(2p+1)$-current in $N$ and $H_{2p+1}(N,\mathbb{Z})=H_{n-(2p+1)}(N,\mathbb{Z})=0$, where $H_{2p+1}(N,\mathbb{Z})$ is the $(2p+1)th$ homology group of $N$, where $2p+1+q=dimN=n$ and $2m-1-dimN=q$.
    \end{thm}
    \begin{proof}
    First of all for a contact CR- warped product submanifold, we have the following
    \begin{eqnarray}
    g(h(\phi X,Z),\phi W)=X(lnf)g(Z,W),\label{eq:4.18}
    \end{eqnarray}
    for $X\in \Gamma(D)$ and $Z\in \Gamma(D^{\perp})$, See:\cite{HM1}. Now using (\ref{eq:4.18}), we also have
    \begin{eqnarray}
    g(h'(X,Z),\phi W)=-\phi X(lnf)g(Z,W).\label{eq:4.19}
    \end{eqnarray}
    Hence for $dim(D^{\perp})=dim(TM^{\perp})$, we derive
    \begin{eqnarray}
    \sum_{i=1}^{2p+1}\sum_{\alpha=1}^{q}g(h'(e_i,e_{\alpha}),h'(e_i,e_{\alpha}))&=& \sum_{i=1}^{2p}\sum_{\alpha=1}^{q}g(h'(e_i,e_{\alpha}),h'(e_i,e_{\alpha}))\nonumber\\
    &&+g(h'(e_{\alpha},\xi),h'(e_{\alpha},\xi)).\nonumber
    \end{eqnarray}
    Since $h'(Y,\xi)=-\phi Y$ (see:\cite[Page 102, Lemma 1.2]{Bejancu}) for $Y\in \Gamma(D^{\perp})$, we obtain
    \begin{equation}
    \sum_{i=1}^{2p+1}\sum_{\alpha=1}^{q}g(h'(e_i,e_{\alpha}),h'(e_i,e_{\alpha}))=\sum_{i=1}^{2p}\sum_{\alpha=1}^{q}g(h'(e_i,e_{\alpha}),h'(e_i,e_{\alpha}))+q. \label{eq:Bay4}
    \end{equation}
    Using (\ref{eq:4.18}) and (\ref{eq:4.19}), we derive
     \begin{eqnarray}
   \sum_{i=1}^{2p+1}\sum_{\alpha=1}^{q}g(h'(e_i,e_{\alpha}),h'(e_i,e_{\alpha}))&=&\sum_{i=1}^{p}\sum_{\alpha,\gamma=1}^{q}(-\phi e_i(lnf))^2g(e_\alpha,e_\gamma)^2\nonumber\\
   &&+e_i(lnf)^2g(e_\alpha,e_\gamma)^2\nonumber\\
   &=&\|\nabla lnf\|^2q. \label{eq:4.20}
    \end{eqnarray}
    Since
 $$
   \sum_{i=1}^{2p+1}\sum_{\alpha=1}^{q}g(Ae_i,e_i)g(Ae_{\alpha},e_{\alpha})=(2p+1-\beta)q,$$
and
   $$\sum_{i=1}^{2p+1}\sum_{\alpha=1}^{q}g(Ae_{\alpha},e_i)g(Ae_{\alpha},e_i)=0,$$
 we derive
    using (\ref{eq:Bay4}), (\ref{eq:3.17}) and (\ref{eq:4.20}), we have
    \begin{eqnarray*}
   && \sum_{i=1}^{2p+1}\sum_{\alpha=1}^{q}2\parallel h'(e_i,e_\alpha),h'(e_i,e_\alpha)\parallel^2-g(h'(e_\alpha,e_\alpha),h'(e_i,e_i))\\
    &&=(2p+1-\beta)q+ \sum_{i=1}^{2p+1}\sum_{\alpha=1}^{q}g(R'(e_i,e_\alpha)e_i,e_\alpha)+\|\nabla lnf\|^2q+q.
    \end{eqnarray*}
    On the other hand since $N$ is a contact CR-submanifold of a Sasakian manifold and it is a warped product submanifold, we have (see \cite{ChenB})
    r
   Hence we get
    \begin{eqnarray*}
    \sum_{i=1}^{2p+1}\sum_{\alpha=1}^{q}g(R'(e_i,e_\alpha)e_i,e_\alpha)&=&\frac{q}{f}\bigtriangleup f.
    \end{eqnarray*}
    Thus we arrive at
    \begin{eqnarray*}
    \sum_{i=1}^{2p+1}\sum_{\alpha=1}^{q}(2\|h'(e_i,e_\alpha)\|^2&-&g(h'(e_\alpha,e_\alpha),h'(e_i,e_i)))=(2p+1-\beta)q\nonumber\\
    &&+\parallel\nabla lnf\parallel^2q+\frac{q}{f}\bigtriangleup f+q.\nonumber
    \end{eqnarray*}
    Therefore, it follows that
    $$\sum_{i=1}^{2p+1}\sum_{\alpha=1}^{q}(2\|h'(e_i,e_\alpha)\|^2-g(h'(e_\alpha,e_\alpha),h'(e_i,e_i)))<2(2p+1)(2p+1-n)\lambda(x)^2$$
    for $x\in N$ if and only if the following inequality is satisfied
      \begin{eqnarray*}
    \bigtriangleup f&<&-\frac{f}{q}\{(2p+1-\beta)q+q\|\nabla lnf\|^2+q-2(2p+1)(2p+1-n)\lambda (x)^2\}
    \end{eqnarray*}
    which completes proof due to $\beta=\eta(A\xi)+1$.
   \end{proof}
   Thus we also have the following result.
  \begin{cor} Let $N^n$ be a compact contact CR-warped product submanifold $M=M_1\times _fM_2$ of a Sasakian hypersurface of $\mathbb{C}^{m}$ such that $M_1$ is a $(2p+1)$-dimensional invariant submanifold tangent to $\xi$ and $M_2$ a $q$-dimensional $\mathcal{C}$- totally real submanifold of $M$. If
    \begin{eqnarray*}
    \bigtriangleup f&<&\frac{f}{q}\{2(2p+1)(2p+1-n)\lambda (x)^2-q(2p-\eta(A\xi))+\|\nabla lnf\|^2+1)\}
    \end{eqnarray*}
    is satisfied for every $x\in N$, where $\lambda (x)^2$ is the maximum at $x$ of the square of principal curvature of $M^{2m-1}$, then $N^n$ is
homeomorphic to a sphere when $n \neq 3$ and  $N^n$ is homotopic to a sphere
when $n = 3$.
\end{cor}
\section{Homology of Contact CR-submanifolds of a Sasakian hypersurface}
In this section, we will consider an arbitrary contact CR-submanifold of a Sasakian hypersurface $M^{2m-1 }$ immersed in $\mathbb{C}^{m}$. By definition of contact CR-submanifold, we use two differentiable distributions in our study. Therefore, we will recall some notions related to distributions which are introduced recently in \cite{R1} and \cite{RZ}. Let $D$ and $D^{\perp}$ be two distributions on a submanifold $M$ of a Riemannian manifold $\bar{M}$. Let $X^T$
be the $D-$ component of $X\in \chi(M)$ (resp., $X^{\perp}$ the $D^{\perp}$-component of $X$).  Then the second fundamental form of $D$ and $D^{\perp}$, respectively, are defined as
 $$h^{\perp}(u,v)=\frac{1}{2}(\nabla_uv+\nabla_vu)^{T}, h^{T}(X,Y)=\frac{1}{2}(\nabla_XY+\nabla_YX)^\perp$$
 for $u,v\in \Gamma(D^{\perp})$ and $X, Y \in \Gamma(D)$.  The mean curvature of $D^{\perp}$ and $D$  are defined as\\
 $$H^T=trace|_gh^T, H^{\perp}=trace|_gh^{\perp}.$$\\
  We also define the following maps
 $$T^{\perp}(u,v)=\frac{1}{2}[u,v]^{T}, T^{T}(X,Y)=\frac{1}{2}[X,Y]^{\perp}$$
  where $T^{\perp}$ and $T^T$ denote the components of $D^{\perp}$ and $D$, respectively. Finally mixed scalar curvature of submanifold of $M$ is given
  $$S_{mix}=\sum_{a,i}g(R(E_\alpha,\varepsilon_i)E_\alpha,\varepsilon_i)$$
    We note that this notion is not depend on the order of distribution on the choice of a local frame.
  By direct computation, we have \cite{R1}
  \begin{eqnarray}
  S_{mix}=\|H^{\perp}\|^2-\|h^{\perp}\|^2+\|T^{\perp}\|^2+\|H^T\|^2-\|h^T\|^2+div(H^T+H^{\perp}).\label{eq:3.16}
  \end{eqnarray}

In this section, we are going to investigate certain conditions for a contact CR-submanifold of a Sasakian hypersurface such that its $2p$-th homology group is zero and as a result of this we obtain a sphere theorem for contact CR-submanifold.
\begin{thm}\label{theo1}
    Let $N$ be an $n$-dimensional compact contact CR-submanifold of a Sasakian hypersurface $M^{2m-1}$ of $\mathbb{C}^{m}$. If $p=dim(D)$ and $q=dim(D^{\perp})$, then
    \begin{eqnarray*}
    S_{mix}&<&2p(\eta(A\xi)-q)+\sum_{i=1}^{p}\sum_{s=1}^{q}-\|(\nabla'_{e_s^*}\phi)e_i\|^2+2g((\nabla'_{e_s^*}\phi)e_i,\phi h'(e_s^*,e_i))\\
    &-&2\|h'(e_i,e_s^*)\|^2+4p(2p-n)\lambda(x)^2
    \end{eqnarray*}
    is satisfied for every $x\in N$ and orthonormal basis $\{e_i, e_s^*\}$ of $D\oplus D^{\perp}$, $i=1,...,2p$, $s=1,...,q$. Then there is no integral $2p$-current in $N$ and $H_{2p}(N,\mathbb{Z})=H_{n-2p}(N,\mathbb{Z})=0$, where $H_p(N,\mathbb{Z})$ is the $p$-th singular homology group with integer coefficients and $h'$ is the second fundamental form of $N$ in $M$.
    \end{thm}
    \begin{proof}
  From (\ref{eq:2.15}) and (\ref{eq:2.6}), we have
  \begin{eqnarray}
 g(h'(X,W),h'(Y,Z))&=&g(AX,Z)g(AY,W)-g(AY,Z)g(AX,W)\nonumber\\
 &&+g(R'(X,Y)Z,W)+g(h'(Y,W),h'(X,Z)).\nonumber
  \end{eqnarray}
  for any vector fields $X, Y, Z$ and $W$ on $N$. We take $X=Z=e_i$, $Y=W=e_{\alpha}$, $\alpha=1,...,q,q+1$
  \begin{eqnarray}
  \sum_{i=1}^{2p}\sum_{\alpha=1}^{q+1}g(h'(e_i,e_{\alpha}),h'(e_{\alpha},e_i))&=&\sum_{i=1}^{2p}\sum_{\alpha=1}^{q+1}(g(Ae_i,e_i)g(Ae_{\alpha},e_{\alpha})-g(Ae_{\alpha},e_i)g(Ae_i,e_{\alpha})\nonumber\\
  &&+g(R'(e_i,e_{\alpha})e_i,e_{\alpha})+g(h'(e_{\alpha},e_{\alpha}),h'(e_i,e_i))).\nonumber
  \end{eqnarray}
Since
 $$
   \sum_{i=1}^{2p}\sum_{\alpha=1}^{q+1}g(Ae_i,e_i)g(Ae_{\alpha},e_{\alpha})=-2p((-q+1)+\beta)=-2p(\beta-(q+1)),$$
and
   $$\sum_{i=1}^{2p}\sum_{\alpha=1}^{q+1}g(Ae_{\alpha},e_i)g(Ae_{\alpha},e_i)=0,$$
 we derive
  \begin{eqnarray*}
  \sum_{i=1}^{2p}\sum_{\alpha=1}^{q+1}g(h'(e_i,e_{\alpha}),h'(e_i,e_{\alpha}))=&-&2p(\beta-(q+1))+\sum_{i=1}^{2p}\sum_{\alpha=1}^{q+1}\{g(R'(e_i,e_{\alpha})e_i,e_{\alpha})\nonumber\\
  &&+g(h'(e_{\alpha},e_{\alpha}),h'(e_i,e_i))\}.
  \end{eqnarray*}
Thus we get
  \begin{eqnarray}
  \sum_{i=1}^{2p}\sum_{\alpha=1}^{q+1}\{2\parallel h'(e_i,e_{\alpha})\parallel^2 &-&g(h'(e_{\alpha},e_{\alpha}),h'(e_i,e_i))\}=-2p(\beta-(q+1))+S_{mix}\nonumber\\
  &&+\sum_{i=1}^{2p}\sum_{\alpha=1}^{q+1}g(h'(e_i,e_{\alpha}),h'(e_i,e_{\alpha})).\label{eq:3.17}
  \end{eqnarray}
Since $h'(e_i,\xi)=0$, we obtain
   \begin{eqnarray}
   \sum_{i,\alpha}g(h'(e_i,e_{\alpha}),h'(e_i,e_{\alpha}))=\sum_{i,s}g(h'(e_i,e_s^*),h'(e_i,e_s^*))\label{eq:Bay1}
    \end{eqnarray}
    Putting (\ref{eq:Bay1}) in (\ref{eq:3.17}), we get
    \begin{eqnarray}
  \sum_{i=1}^{2p}\sum_{\alpha=1}^{q+1}\{2\parallel h'(e_i,e_{\alpha})\parallel^2 &-&g(h'(e_{\alpha},e_{\alpha}),h'(e_i,e_i))\}=-2p(\beta-(q+1))+S_{mix}\nonumber\\
  &&+\sum_{i=1}^{2p}\sum_{\alpha=1}^{q}g(h'(e_i,e^*_{\alpha}),h'(e_i,e^*_{\alpha})).\nonumber
  \end{eqnarray}
    On the other hand, by direct calculation we have
    \begin{eqnarray}
    h'(\phi X,Z)=-(\nabla'_{Z}\phi)X+\phi h'(Z,X),\nonumber
    \end{eqnarray}
    where $\nabla'$ is the induced Levi-Civita connection on $M$. Thus we have
    \begin{eqnarray}
    g(h'(\phi X,Z),h'(\phi X,Z))&=&\|(\nabla'_{Z}\phi)X)\|^2-2g((\nabla'_{Z}\phi)X,h'(Z,X))\nonumber\\
    &&+g(\phi h'(Z,X),\phi h'(Z,X)).\nonumber
    \end{eqnarray}
    for $X\in \Gamma(D)$ and $Z\in \Gamma(D^{\perp})$. Since $ \eta( h'(Z,X)=0$, we obtain
    \begin{eqnarray}
    g(\phi h'(Z,X),\phi h'(Z,X))= g( h'(Z,X), h'(Z,X)))=g( h'(Z,X), h'(Z,X)).\nonumber
    \end{eqnarray}
   Putting above equality in (\ref{eq:Bay1}), we have
    \begin{eqnarray*}
    \sum_{i=1}^{2p}\sum_{s=1}^{q}g(h'(e_i,e_s^*),h'(e_i,e_s^*))&=&\sum_{i=1}^{p}\sum_{s=1}^{q}(\|(\nabla_{e_s^*}\phi)e_i\|^2-2g((\nabla'_{e_s^*}\phi)e_i,\phi h'(e_s^*,e_i))\\
    &&+2g(h'(e_i,e_s^*),h'(e_i,e_s^*))).
    \end{eqnarray*}
    Thus we arrive at
    \begin{eqnarray}
    \sum_{i=1}^{2p}\sum_{\alpha=1}^{q+1}(2\|h'(e_i,e_{\alpha})\|^2&-&g(h'(e_{\alpha}e_{\alpha}),h'(e_i,e_i)))=-2p(\beta-(q+1))+S_{mix}\nonumber\\
    &+&\sum_{i=1}^{p}\sum_{s=1}^{q}(\|(\nabla'_{e_s^*}\phi)e_i\|^2-2g((\nabla'_{e_s^*}\phi)e_i,\phi h'(e_s^*,e_i))\nonumber\\
    &&+2\|h'(e_i,e_s^*)\|^2)\label{eq:Bay3}.
    \end{eqnarray}
    From above equality and Theorem \ref{Cheng-Shiohama}, the proof is complete due to $\beta=\eta(A\xi)+1$.
    \end{proof}
    From Theorem \ref{theo1}, (\ref{eq:Bay3}) and (\ref{eq:3.16}), we have the following result.

    \begin{cor}
    Let $N$ be a $n$-dimensional compact CR-submanifolds of $M^{2m-1}$ which is a $(2m-1)$-dimensional Sasakian hypersurface of $\mathbb{C}^{m}$. If $p=dim(D)$ and $q=dim(D^{\perp})$, then
    \begin{eqnarray*}
    \|H^{\perp}\|^2+\|T^{\perp}\|^2+\|H^T\|^2+div(H^T+H^{\perp})+\sum_{i=1}^{p}\sum_{s=1}^{q}(\|(\nabla'_{e_s^*}\phi)e_i\|^2
    +2\|h'(e_i,e_s^*\|^2)\\
    <2p(\eta(A\xi)-q))+4p(2p-n)\lambda^2+\sum_{i=1}^{p}\sum_{s=1}^{q}2g((\nabla'_{e_s^*}\phi)e_i,\phi h'(e_s^*,e_i)).
    \end{eqnarray*}
    is satisfied for every $x\in N$ and orthonormal basis $\{e_i,...,e_s^*\}$ of $D\oplus D^{\perp}$, $i=1,...,2p$, $s=1,...,q$. Then there is no integral $2p$-current in $N$ and $H_{2p}(N,\mathbb{Z})=H_{n-2p}(N,\mathbb{Z})=0$, where $H_p(N,\mathbb{Z})$ is the $p$-th singular homology group with integer coefficients and $h'$ is the second fundamental form of $N$ in $M$.
    \end{cor}
    \begin{proof}
    By using (\ref{eq:3.16}), (\ref{eq:Bay3}) becomes
    \begin{eqnarray*}
    \sum_{i=1}^{2p}\sum_{\alpha=1}^{q+1}(2\|h'(e_i,e_{\alpha}\|^2&-&g(h'(e_{\alpha}e_{\alpha}),h'(e_i,e_i)))=-2p(\beta-(q+1))\\
    &&+\|H\|^2-\|h\|^2+\|T\|^2+\|H^{\perp}\|^2-\|h^{\perp}\|^2\\
    &&+div(H+H^{\perp})+\sum_{i=1}^{p}\sum_{s=1}^{q}(\|(\nabla'_{e_s^*}\phi)e_i\|^2\\
    &&-2g((\nabla'_{e_s^*}\phi)e_i,\phi h'(e_s^*,e_i))+2\|h'(e_i,e_s^*)\|^2).
    \end{eqnarray*}
    By using $\beta=\eta(A\xi)+1$ we obtain the assertion.
    \end{proof}
    We also have the following result.
    \begin{cor}
    Let $N$ be an $n$-dimensional compact contact CR-submanifold of a Sasakian hypersurface $M^{2m-1}$ of $\mathbb{C}^{m}$. If $p=dim(D)$ and $q=dim(D^{\perp})$, then
    \begin{eqnarray*}
    S_{mix}&<&2p(\eta(A\xi)-q))+\sum_{i=1}^{p}\sum_{s=1}^{q}(-\|(\nabla'_{e_s^*}\phi)e_i\|^2+2g((\nabla'_{e_s^*}\phi)e_i,\phi h'(e_s^*,e_i))\\
    &-&2\|h'(e_i,e_s^*)\|^2)+4p(2p-n)\lambda(x)^2
    \end{eqnarray*}
    is satisfied for every $x\in N$ and orthonormal basis $\{e_i,...,e_s^*\}$ of $D\oplus D^{\perp}$, $i=1,...,2p$, $s=1,...,q$,  then $N^n$ is
homeomorphic to a sphere when $n \neq 3$ and  $N^n$ is homotopic to a sphere
when $n = 3$,  where $h'$ is the second fundamental form of $N$ in $M$.
\end{cor}

\section{Concluding Remarks}
In this article, we have taken a step in understanding the homology groups of the contact CR-submanifold of a Sasakian hypersurface as the first step in understanding the homology groups of a general contact CR-submanifold of an arbitrary almost contact manifold. Our study will stimulate to obtain more results  for more general situations.\\

\noindent{\bf Acknowledgment.} This work was supported by
Research Fund of the Ege University with Project Number:20776.

\end{document}